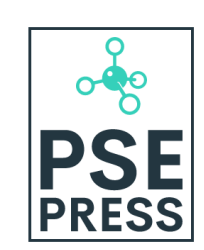

# Enhancing Interpretability of Stochastic Programming Solutions: A Multiparametric Approach

**Parth Brahmbhatt[a] and Styliani Avraamidou[a]***

[a] Department of Chemical and Biological Engineering, University of Wisconsin-Madison, Madison, WI, USA
* Corresponding Author: avraamidou@wisc.edu

## ABSTRACT

Stochastic programming (SP) is a powerful framework for decision-making under uncertainty, but its practical adoption in industry is often hindered by the difficulty in understanding the causal relationships that drive optimal solutions. In the two-stage SP, strategic first-stage decisions are coupled with operational second-stage recourse decisions. When the number of scenarios under consideration is large, understanding the direct link between the uncertainty realization and optimal recourse strategy becomes computationally and cognitively demanding. Common approaches to improve interpretability include trained classification trees or scenario reduction, replacing the large scenario set with a representative subset. This is often achieved through post-hoc clustering (e.g., k-means) based on uncertainty realizations or optimal recourse decisions. While useful, these methods only provide a statistical approximation of the solution space and may fail to reveal the underlying structural properties of the recourse problem that drive optimal first-stage decisions. This work introduces a novel, deterministic approach to explainability using multiparametric programming (mp) within a Benders decomposition framework. We reformulate the recourse subproblem as a multiparametric linear program, generating an explicit map of Critical Regions (CRs), which are polyhedral partitions of the uncertainty space. This allows us to cluster scenarios analytically rather than statistically. We demonstrate this methodology on a supply chain planning under demand uncertainty. Our results show that 100 stochastic scenarios map to exactly six critical region clusters. This mapping allows us to explain optimal capacity planning decisions as a precise trade-off between specific operational modes, providing a fully transparent interpretation of the stochastic solution.

## INTRODUCTION

Optimization under uncertainty is central to Process Systems Engineering (PSE), particularly in applications such as supply chain management, energy systems, and capacity planning. Two-stage stochastic programming (SP) [1] is the standard approach for these problems, distinguishing between "here-and-now" decisions (made before uncertainty reveals itself) and "wait-and-see" recourse decisions (made after the uncertainty is realized).

Despite the mathematical maturity of SP, the resulting solutions often function as "black boxes" for decision-makers [2]. A typical large-scale SP solution provides a single optimal first-stage vector $x^*$ that minimizes expected cost across thousands of scenarios. However, it rarely explains why this specific $x^*$ is robust. For instance, in a capacity expansion problem, a manager might ask: "Why is the optimal capacity for Plant B set to 26 units rather than 30?" The answer lies buried in the complex interactions of thousands of second-stage variables.

Recent developments in Explainable AI [3] have attempted to bridge this gap. In the context of optimization, explainability usually relies on sensitivity analysis, which is inherently local, or post-hoc analysis using machine learning surrogates like classification trees [4].

A common approach to reduce complexity is scenario reduction [5], where a large set of scenarios is clustered into a smaller representative set. This is often achieved through post-hoc clustering [6] techniques (e.g., k-means). Rathi et al. [2] recently proposed

"recourse-based clustering, " where scenarios are grouped not by their input uncertain parameters (e.g., demand), but by the similarity of their optimal recourse decisions ($y^*$). While the method is improved compared to input-based clustering, it remains a statistical approximation, requiring the user to arbitrarily select the number of clusters ($k$), potentially obscuring the exact boundaries where the optimal decision strategy changes.

In this work, we propose a rigorous, analytical alternative to statistical clustering. We leverage the synergy between Benders Decomposition [7] and Multiparametric Programming (mp) [8]. By treating the first-stage variables and uncertainty as parameters in the second-stage subproblem, we can solve the subproblem explicitly [9]. The result is a set of Critical Regions (CRs). Each CR represents a partition of the parameter space where the optimal operating strategy (the set of active constraints) is identical. This allows us to:

- **Derive Explicit Policy Laws:** Within each cluster, the recourse decision is an explicit affine function of the uncertainty.
- **Cluster Scenarios Analytically:** Scenarios are grouped if they fall into the same CR.
- **Quantify Trade-offs:** We can explain first-stage decisions by summing the probability-weighted contributions of specific CRs.

# PRELIMINARIES

## Two-Stage Stochastic Programming

We consider the standard two-stage SP formulation, which can be a Linear Program (LP) or a Mixed Integer Linear Program (MILP):

$$\begin{aligned} &\min_{\mathbf{x}} c^{\top}\mathbf{x} + \mathbb{E}_{\xi}[\mathcal{Q}(\mathbf{x},\xi)] \\ &\text{s.t.}\, A\mathbf{x} \le b;\ \mathbf{x} \ge 0 \end{aligned} \tag{1}$$

where $\mathbf{x}$ represents the first-stage decisions (e.g., capacity expansion). The vector $c$ represents the cost coefficients, while $A$ and $b$ define the first-stage constraints. The recourse function $\mathcal{Q}(\mathbf{x},\xi)$ is defined as the optimal value of the second-stage problem for a specific uncertainty realization $\xi$:

$$\begin{aligned} &\mathcal{Q}(\mathbf{x},\xi) = \min_{\mathbf{y}} q^{\top}\mathbf{y} \\ &\text{s.t.}\, W\mathbf{y} \le h(\xi) - T\mathbf{x};\ \mathbf{y} \ge 0 \end{aligned} \tag{2}$$

Here, $\mathbf{y}$ represents operational decisions (e.g., production rates, flows), and $\xi$ captures uncertainties such as demand or prices. The $q$ is the objective coefficient vector, and $W$ (the recourse matrix), $T$ (the technology matrix), and $h$ (right-hand side uncertainty) define the recourse constraints.

## Benders Decomposition

For large-scale problems, the two-stage SP is typically solved using a decomposition approach, such as Benders Decomposition [7], which partitions the problem into a Master problem and a set of sub-problems. The Master problem approximates the recourse cost using cuts [7]:

$$\begin{aligned} &\min_{\mathbf{x},\eta} c^{\top}\mathbf{x} + \sum_{\omega\in\Omega} p_{\omega}\eta_{\omega} \\ &\text{s.t.}\, \eta_{\omega} \ge \lambda_{\omega}^{\top}(h_{\omega} - T\mathbf{x})\ \forall\ \text{cuts} \end{aligned} \tag{3}$$

Where $p_{\omega}$ is probability of scenario $\omega$ and $\lambda$ are dual multipliers from the subproblem. Note that we assume relatively complete recourse holds for the problem. While efficient for computation, Benders decomposition does not inherently generate interpretable insights.

## Multiparametric Programming Formulation

In our recent work [9], we reformulate the subproblem using multiparametric programming. In this work, we are exploring this method to gain interpretable insights.

In the recourse problem, the optimal second-stage decision $\mathbf{y}^*$ depends on two vectors: the fixed first-stage decision $\mathbf{x}$ and the random parameter realization $\xi$. We aggregate these into a single parameter vector $\theta = [\mathbf{x},\xi]^{\top}$ [9]. The recourse problem becomes a multiparametric Linear Program (mp-LP) :

$$\begin{aligned} &z(\theta) = \min_{\mathbf{y}} q^{\top}\mathbf{y} \\ &\text{s.t.}\, W\mathbf{y} = F\theta + d;\ \mathbf{y} \ge 0 \end{aligned} \tag{4}$$

The solution to this mp-LP is not a single scalar value, but a set of critical regions. Each region $CR_i$ is a convex polyhedron defined by a set of linear inequalities:

$$CR_i = \{\theta \in \Theta \mid E_i\theta \le f_i\}$$

where $E_i$ and $f_i$ are a matrix and vector (boundary constraints of $CR_i$), respectively, derived from the optimality conditions (active constraints) of the subproblem. The parameter space $\Theta$ is partitioned into non-overlapping polyhedral regions $CR = \cup_i CR_i$. Within each region $CR_i$, the optimal decision variables $y^*$ and the objective function $z$ are explicit piecewise affine functions of $\theta$:

$$\begin{aligned} &\mathbf{y}^*(\theta) = H_i\theta + k_i, \forall\theta \in CR_i \\ &z^*(\theta) = A_i\theta + B_i, \forall\theta \in CR_i \end{aligned} \tag{5}$$

Further details on this multiparametric-based Benders decomposition framework can be found in Brahmbhatt et al. [9]. This structure provides the "white-box" explainability: every scenario $\xi$ maps to a specific $CR_i$, and that region dictates exactly which constraints are binding (e.g., "Demand constraint is active" vs "Capacity constraint is active").

## METHODOLOGY

Our proposed framework integrates mp-programming into the post-optimality analysis of stochastic programming solutions.

**Step 1: Parametric Solution-** We treat the recourse problem as an mp-LP, considering uncertain parameters (e.g., the first stage decision $\mathbf{x}$ and demand $\xi$) as the varying parameters $\theta$. We solve this using an mp algorithm [11]. This generates the set of all possible operational modes (critical regions) the system can experience under the optimal design.

**Step 2: Solve the Full Stochastic Problem-** We solve the original two-stage SP using either the mp Benders method [9] or a monolithic solver to obtain the optimal first-stage solution $\mathbf{x}^*$.

**Step 3: Map Scenarios to Regions: Clustering-** For a given optimal fixed first-stage decision $\mathbf{x}^*$, we project the full set of sampled scenarios $\Omega$ onto the generated critical regions.

$$\text{Cluster}_i = \{ \omega \in \Omega \mid \xi_\omega \in CR_i \}$$

As shown in Figure 1, this mapping is treated as a point location problem. This process simply iterates through the generated regions to find the specific critical region ($CR_i$) where the boundary constraints are satisfied for a given scenario. Here in Figure 1, the demand of 3 different customers is uncertain in the 2nd stage. So, using parametric space for a fixed first-stage decision, we can cluster all demand scenarios in critical region clusters. Unlike k-means, which clusters based on Euclidean distance, this step clusters scenarios based on structural identity. Two scenarios are in the same cluster if and only if they share the exact same active constraints.

**Step 4: Analytical Interpretation-** We explain the optimal $\mathbf{x}^*$ by analyzing the probability mass accumulated in "stress" regions (regions with high penalties or recourse costs) versus "nominal" regions.

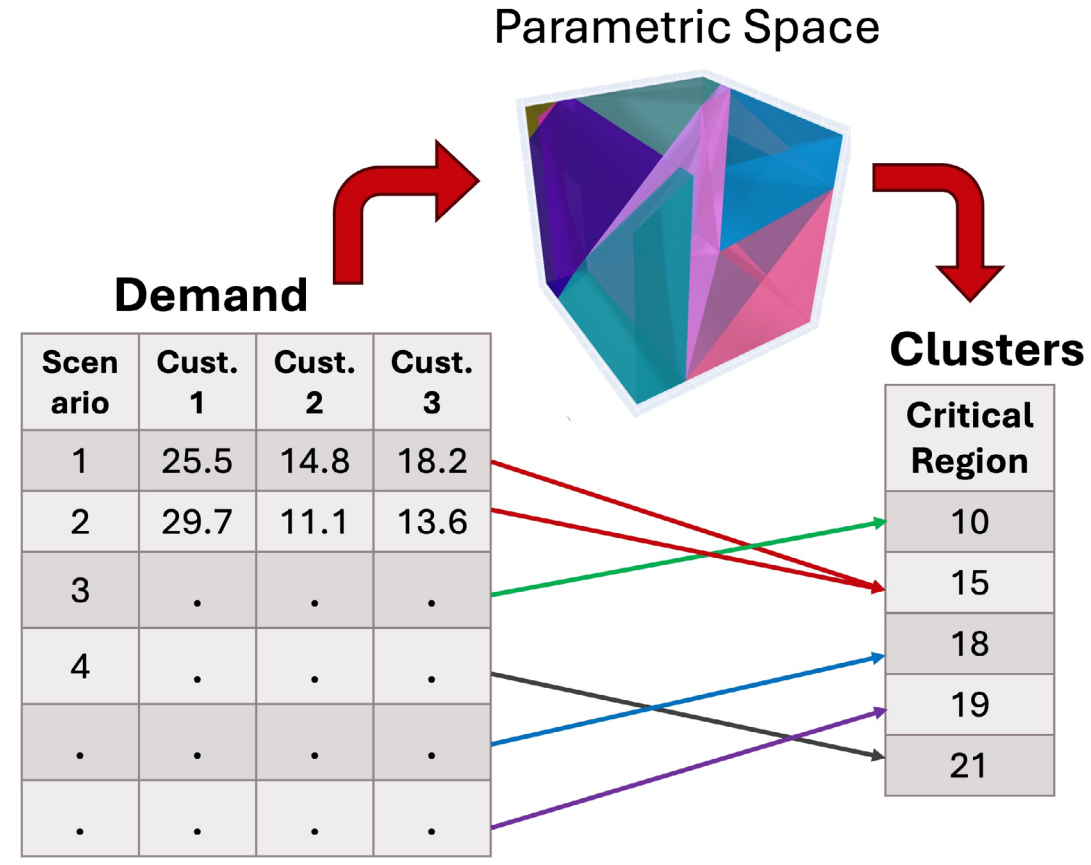

**Figure 1.** Recourse-based clustering of uncertain scenarios using mp critical regions.

## CASE STUDY: SUPPLY CHAIN PLANNING

To demonstrate the efficacy of the proposed multiparametric framework, we utilize a canonical supply chain network design problem adapted from Rathi et al [2]. This problem represents a standard capacity planning challenge where strategic decisions must be made against uncertain customer demands.

### Mathematical Formulation

We consider a network consisting of a set of manufacturing plants $i \in I$ and a set of customers $j \in J$. The problem is formulated as a two-stage stochastic linear program.

**First-Stage Problem:** The decision-maker determines the production quantities $x_i$ for each plant $i$ prior to the realization of uncertain demand. The objective is to minimize the total production cost plus the expected recourse cost. The $C_i^{prod}$ are production costs and $C_i^{max}$ the maximum production capacities of the plants, respectively.

$$\min_x \sum_{i \in I} C_i^{prod} x_i + \mathbb{E}_\xi[\mathcal{Q}(x,\xi)] \qquad (6)$$
$$\text{s.t. } 0 \le x_i \le C_i^{max} \; \forall \, i \in I$$

**Second-Stage (Recourse) Problem:** Given the first-stage production levels $x$ and a demand realization $\xi$ (where $d_j(\xi)$ is the demand for customer $j$), the recourse problem optimizes transportation flows $y_{ij}$ and external procurement $\bar{y}_j$. External procurement represents a penalty mechanism (e.g., outsourcing or lost sales) used when local production cannot satisfy demand.

$$\mathcal{Q}(x,\xi) = \min_{y,\bar{y}} \sum_{i \in I} \sum_{j \in J} C_{ij}^{trans} y_{ij} + \sum_{j \in J} C_j^{ext} \bar{y}_j$$
$$\text{s.t. } \sum_{j \in J} y_{ij} \le x_i \; \forall i \in I \text{ (Supply Capacity)} \qquad (7)$$
$$\sum_{i \in I} y_{ij} + \bar{y}_j \ge d_j(\xi) \; \forall j \in J \text{ (Demand Satisfaction)}$$
$$y_{ij}, \bar{y}_j \ge 0$$

Here, $C_{ij}^{trans}$ is the transportation cost from plant $i$ to customer $j$, and $C_j^{ext}$ is the cost of external procurement, incentivizing satisfying the demand via internal production whenever possible.

### Network Topology and Parameters

We apply this formulation to an instance consisting of two plants ($I = 2$, labeled A and B) and 3 Customers ($J = 3$, labeled 1, 2, and 3) in Figure 2.

**First-Stage Parameters:** The two plants differ in production efficiency. The production costs ($C_A^{prod}, C_B^{prod}$) are \$100/unit and \$150/unit, respectively. The maximum capacity ($C_i^{max}$) is 38.0 units for both plants.

**Second-Stage Parameters:** The unit transportation costs ($C_{ij}^{trans}$) from Plant A to Customers 1, 2, and 3 are \$90, \$67.5, and \$90, respectively. Similarly, the costs from Plant B to Customers 1, 2, and 3 are \$67.5, \$30, and \$67.5, respectively. The external procurement cost ($C^{ext}$) is set at \$500 per unit for each customer. External procurement is set significantly higher than production plus transport to penalize shortages.

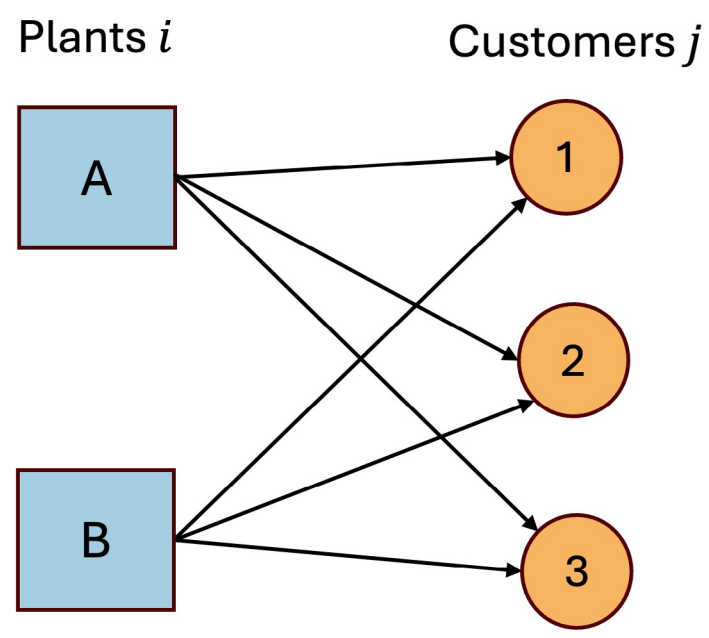

**Figure 2.** The network topology for 2 plants and 3 customers.

**Uncertainty Modeling:** The demand $d_j(\xi)$ for each of the three customers is uncertain. We approximate the uncertainty using a sample average approximation (SAA) approach with $N = 100$ scenarios. The demand for each customer $j$ follows a uniform distribution U(10, 30).

The resulting problem has 2 first-stage variables, 900 second-stage variables (9 per scenario × 100 scenarios), and 700 constraints. While solvable by standard monolithic solvers, the focus of this study is on analyzing the solution structure using the multiparametric decomposition proposed in [9] and summarized in the previous section.

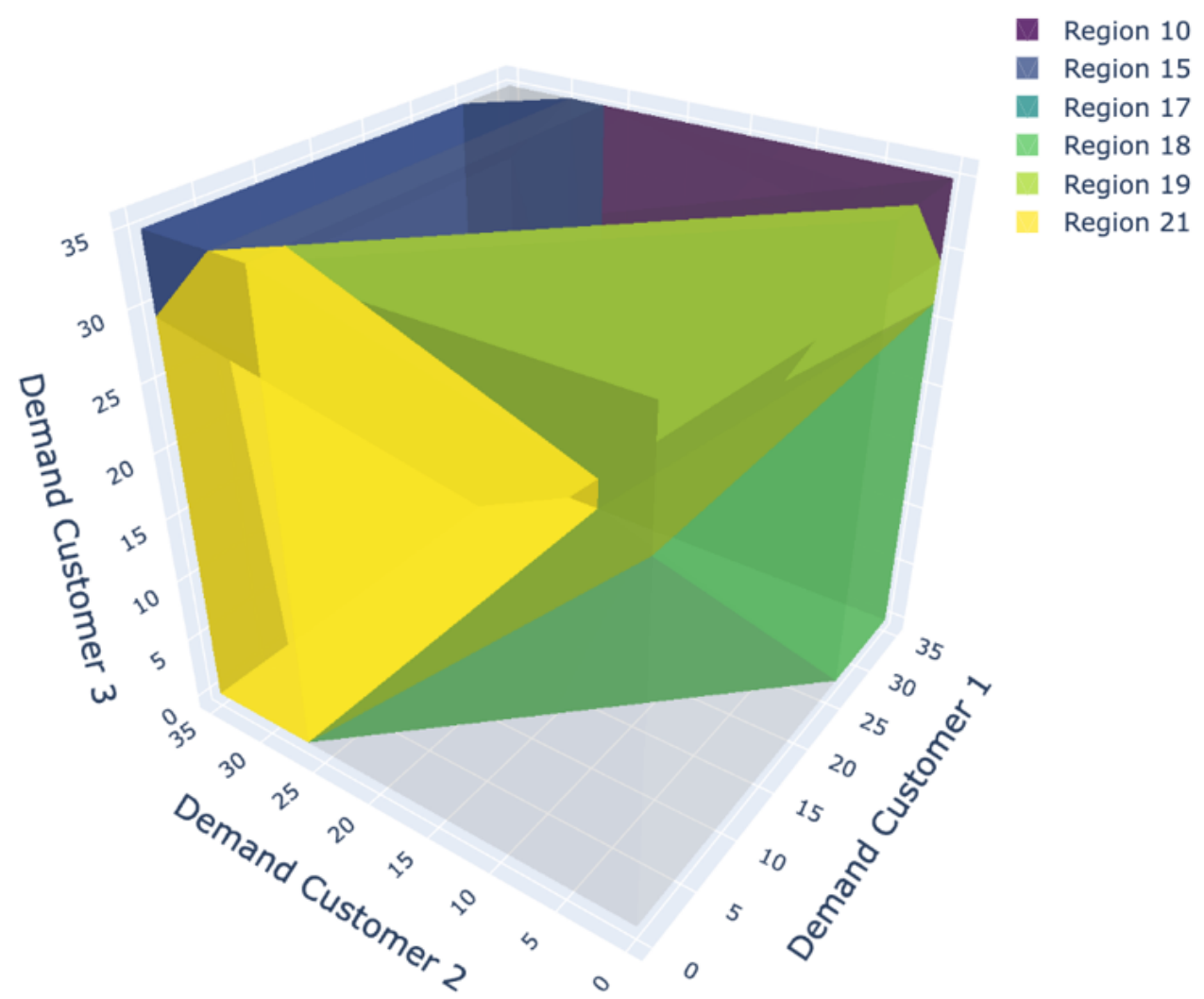

**Figure 3.** Critical regions in demand space given fixed first-stage decisions (each color represents a different critical region).

# RESULTS AND DISCUSSION

The optimal first-stage solution found is Plant A, having a capacity of 38.00 units, and Plant B, having a capacity of 26.41 units, for a total production capacity of 64.41 units.

## Critical Region Clusters Analysis

The multiparametric analysis was implemented in Python using the PPOPT solver [10], and it resulted in 23 critical regions. For the fixed first stage solution ($x^*$), Figure 3 shows the partitioned 3D demand space (Customer 1, 2, 3) with critical regions. Each color in Figure 3 represents a different critical region. While the statistical approach, like k-means, required selecting an arbitrary number of clusters (say k=5 or 6), our analytical approach identified exactly 6 populated critical regions that cover all 100 scenarios. As shown in Figure 1, the proposed clustering method uses the point location problem to find the critical region ID. This is what we call recourse-based clustering, which is shown in Figure 4. Crucially, it demonstrates that while the uncertainty space is continuous, the corresponding optimal recourse strategy is strictly piecewise affine and governed by exactly six distinct critical regions, each defined by a unique set of active constraints in the recourse (sub-problem) decision space.

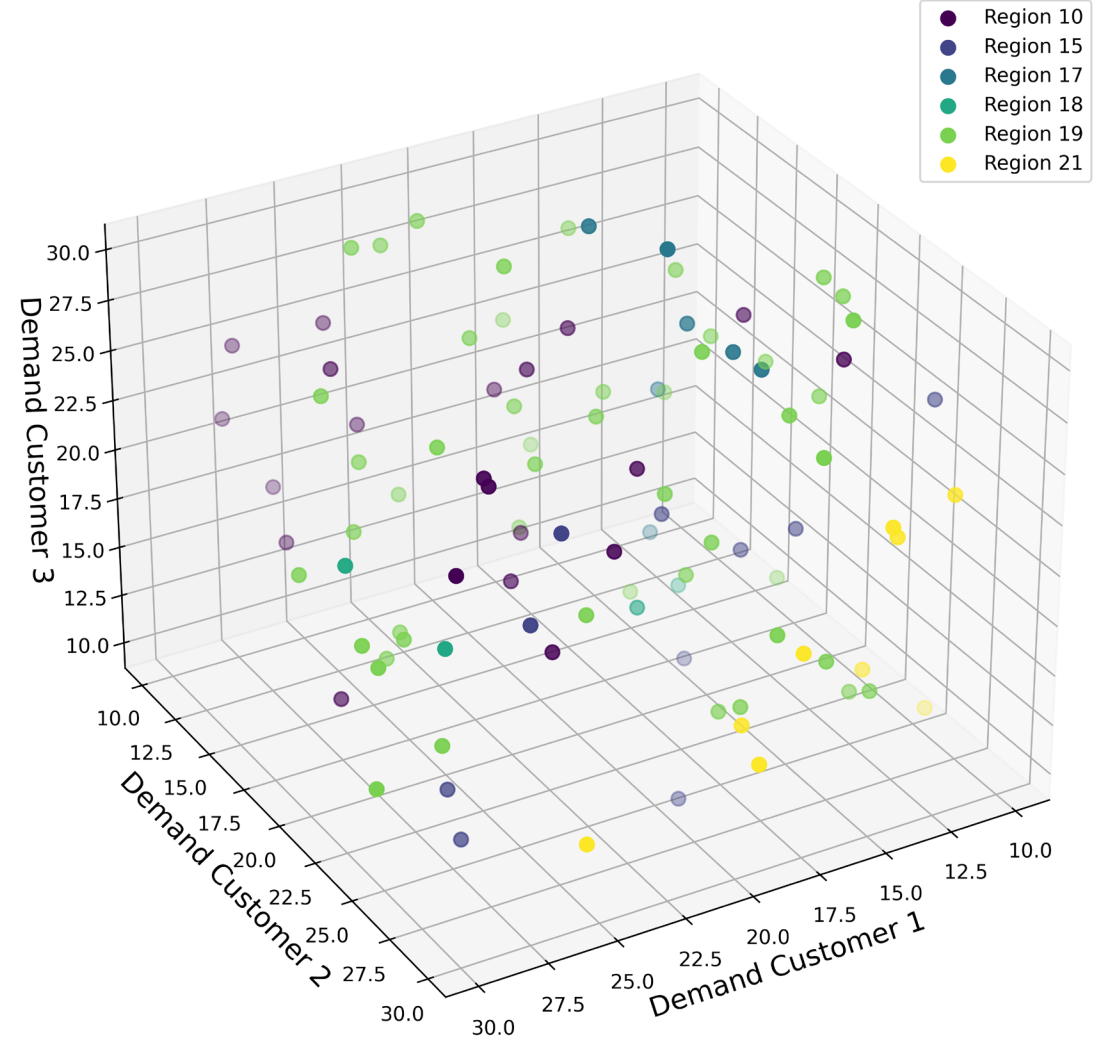

**Figure 4.** 3D Visualization of the 100 demand scenarios mapped to 6 Critical Regions. The recourse clustering is defined by polyhedral boundaries (active constraints) rather than simple distance.

A comparison of our recourse-based clustering with scenario realization-based clustering is presented in Figure 5, which shows the 100 scenarios clustered using the k-means clustering method, indicating that scenarios of the same type were grouped together within the same cluster ID. But based on recourse-based clustering using mp critical region, we can say that scenarios that are similar in their uncertainty realizations may not lead to the same optimal recourse decisions. Figure 6 shows the clusters of recourse decision for all 100 scenarios based on the mp critical regions. Within the same critical region clusters, the scenarios have the same recourse decisions.

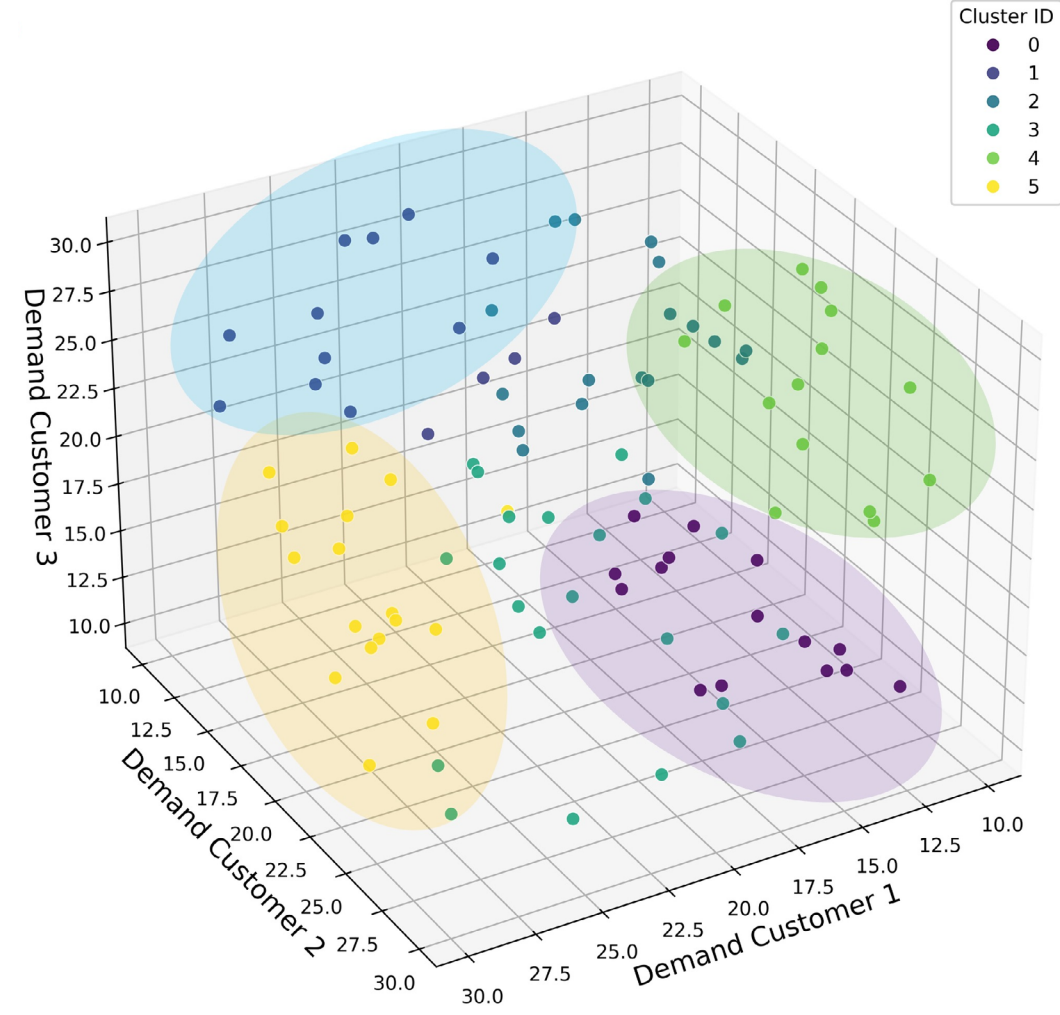

**Figure 5.** 3D Visualization of the 100 demand scenarios mapped to 6 clusters using k-means clustering.

## Operational Modes (Clusters)

Figure 7 shows the average recourse decisions over different critical region clusters. The size of the square shows the production level, which is fixed for all scenarios. The size of the circle shows the average demand in each critical region cluster. As we can see Region 19 (Probability: 0.49) is the "Nominal Operation" mode. It contains nearly half of all scenarios. In this region, total demand is moderate, and the production from A and B is sufficient. No external procurement is required. Region 21 (Probability: 0.09): A transition region where demand is slightly higher, shifting transportation routes, but still manageable without external purchase. Regions 10 & 15 (Probability: 0.21 & 0.10): These are "Failure/Recourse" modes. In these regions, high simultaneous demand on the customer's side exceeds the available inventory (64.41 units). The system hits the active constraint $\sum y \leq x$, forcing the variable $\bar{y}$ (external procurement) to become non-zero. The average amount of this outside procurement is explicitly indicated by the red numbers next to the customer nodes in Figure 7. The probability here is calculated by the number of scenarios that got clustered in each specific region over the total scenarios.

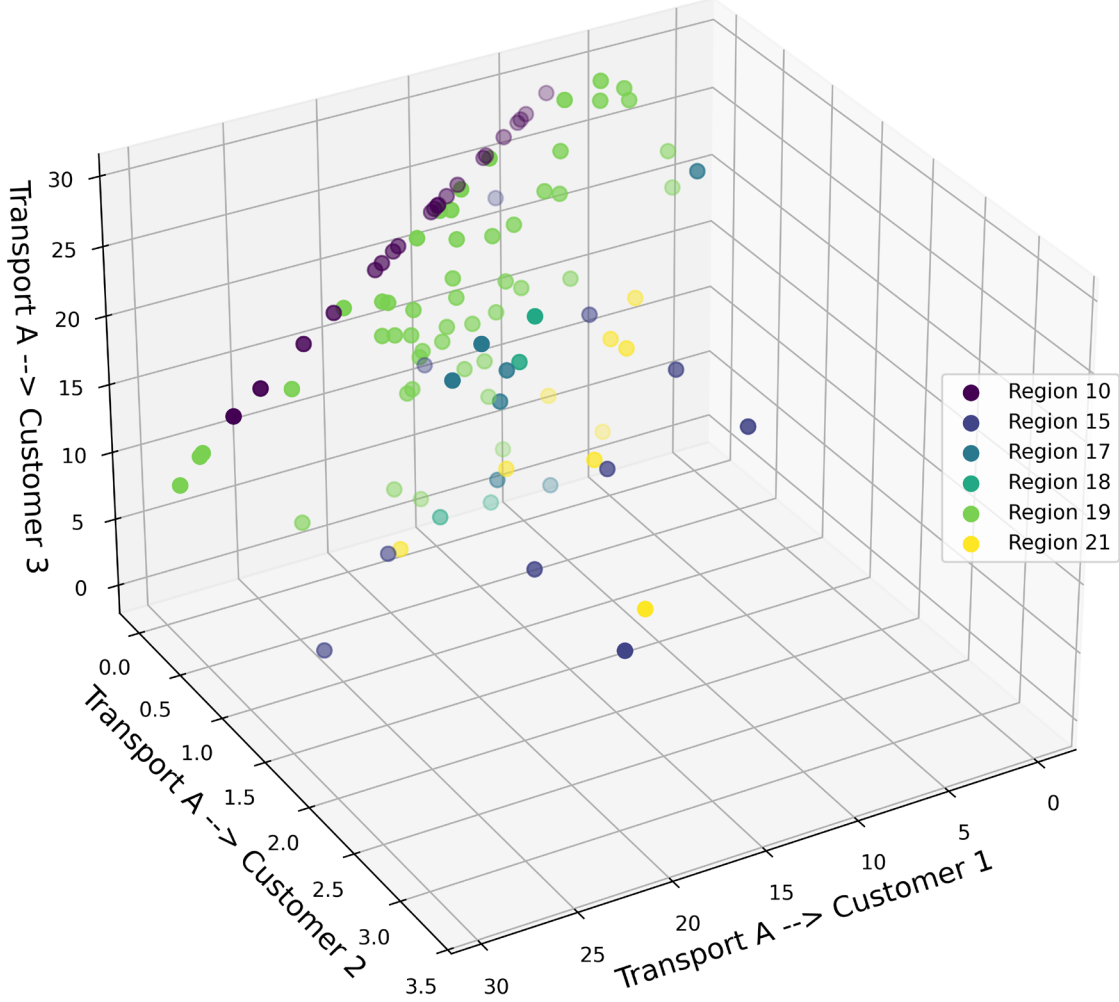

**Figure 6.** The recourse decisions of all 100 scenarios got clustered in 6 mp clusters. The 3 axes represent the transportation of the product from Plant A to customers 1, 2, and 3, respectively.

## Explaining the Capacity Decision

The central explainability question for this SP solution is: Why is the optimal production at Plant B set to 26.41 units instead of its maximum 38.0 units? A human planner might argue: "We are hitting penalty costs in Regions 10 and 15 (31% of the time). Why not increase production at Plant B to avoid this?" Our mp based framework allows us to answer this with a precise cost-benefit calculation. As shown in Figure 8, for the higher demand for all customer cases (Regions 10 and 15), we have very high second-stage costs. But these 2 clusters (Region 10 and 15) have an overall low probability; they make a very low contribution to the overall second-stage cost. Current external procurement cost contribution from Regions 10 and 15 is approximately $2220 (weighted by probability). While increasing Plant B capacity and transportation of product to customers to cover these scenarios would cost roughly $3319 in additional first-stage investment.

It is mathematically cheaper to accept the "failure" (external procurement) in 31% of scenarios falling into Regions 10/15 than to pay for the capacity to cover them. This transforms the result from a simple number into a clear economic decision. The model intentionally incurs penalties in these specific high-demand regions rather than over-investing in capacity.

While the analysis above explains why we do not increase capacity to meet peak demand in Regions 10 and 15, a symmetrical question arises: Why is the optimal capacity for Plant B set to 26.41 units and not lower? The multiparametric framework allows us to answer this by analyzing the "nominal" operating modes, specifically Region 19 and Region 21. These regions have a combined probability of 0.58 and average demands (over all customers) of 56.51 and 60.10 units, respectively. Currently,

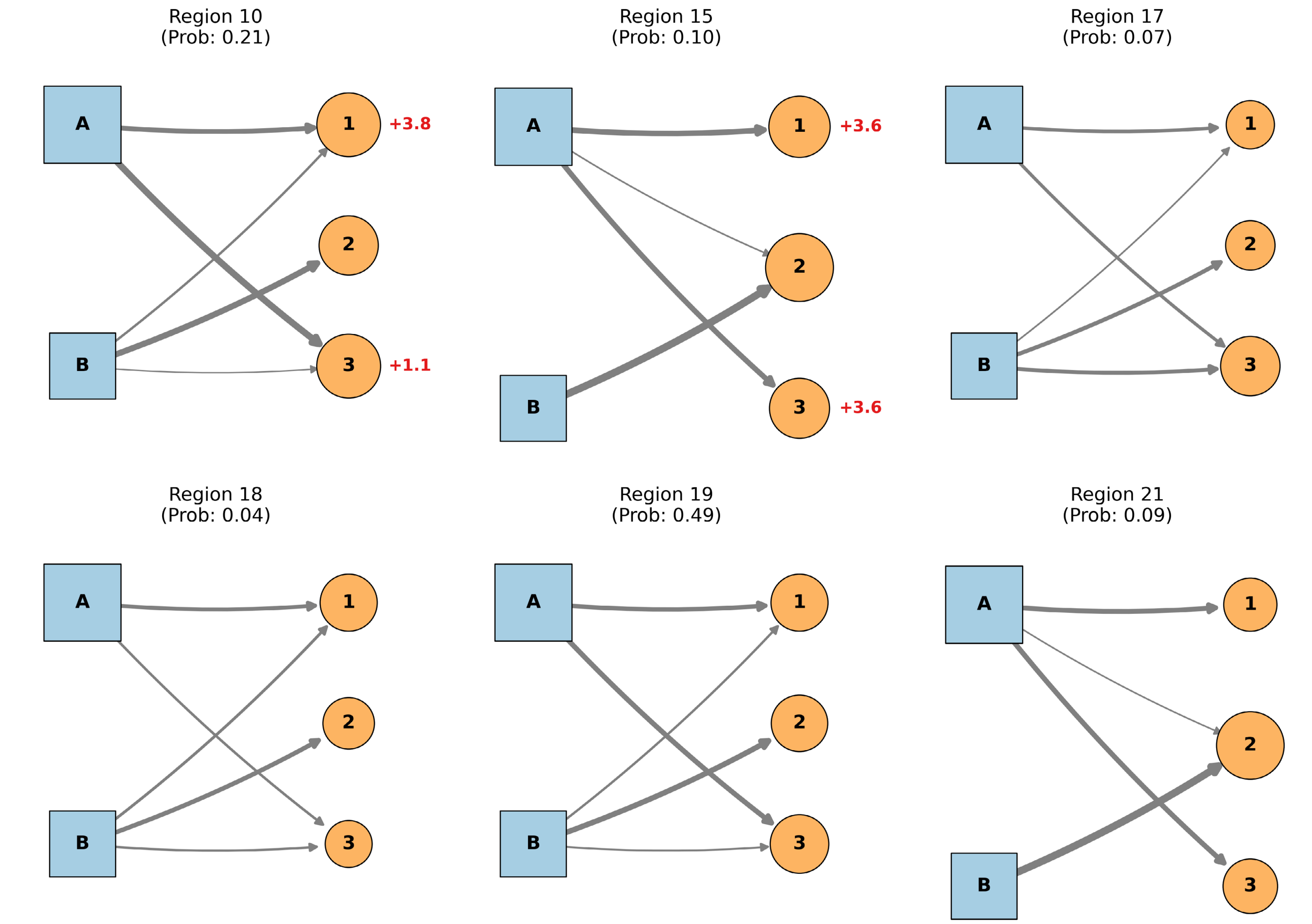

**Figure 7.** Visualization of recourse decisions averaged across the six Critical Region clusters. The square nodes represent manufacturing plants, where the size corresponds to the fixed production capacity. The circular nodes represent customers, with their size proportional to the average demand within that specific cluster. The red numbers adjacent to customer nodes indicate the specific quantity of external procurement required to satisfy

the total system capacity of 64.41 units is sufficient to satisfy these demands without external procurement. However, if we were to further reduce Plant B's capacity, we would force external procurement in these highly probable regions. Reducing capacity by one unit saves the guaranteed investment cost but incurs an expected penalty of \$445.0 in additional external procurement costs across these clusters. This far exceeds the marginal cost of \$182.32 to produce and transport that unit internally from Plant B. Thus, the optimal capacity of 26.41 units represents the precise equilibrium where the marginal savings of capacity reduction are outweighed by the expected penalties in the dominant operating regions.

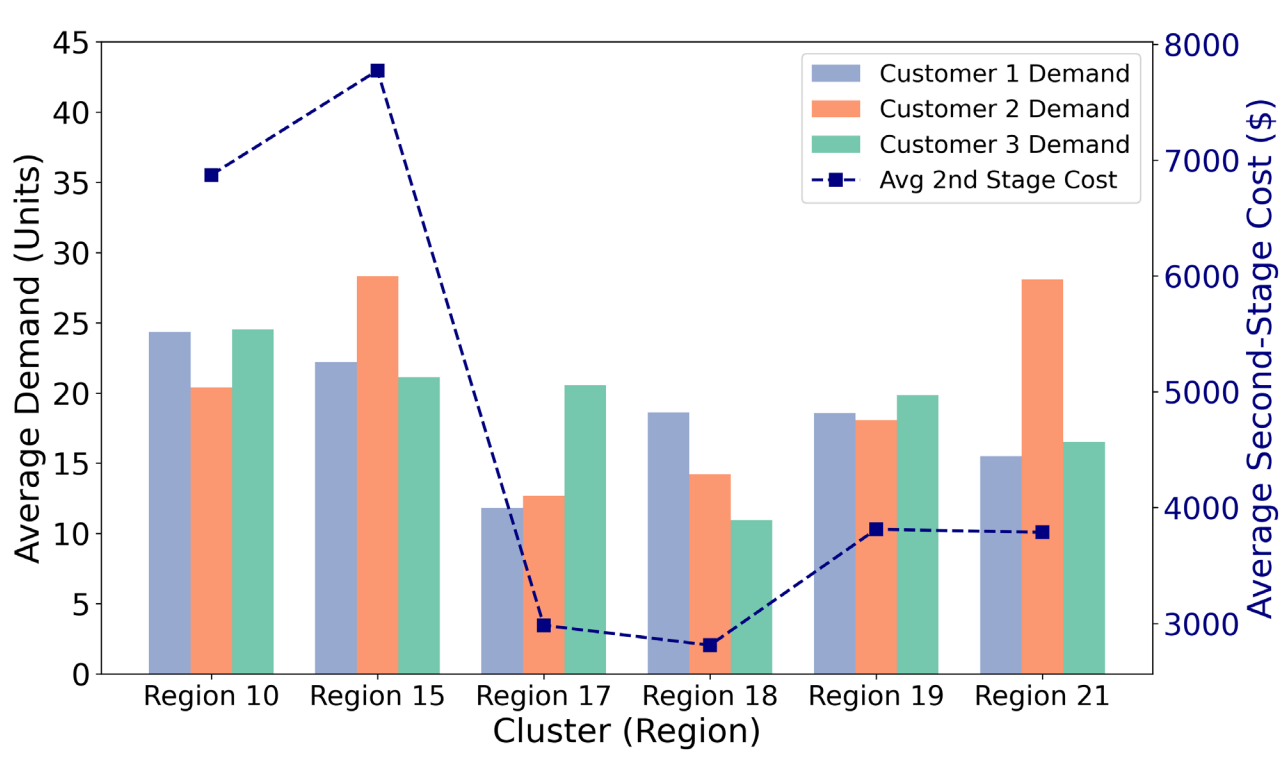

**Figure 8.** Average demand and average second-stage cost in different clusters.

## CONCLUSION

We presented a methodology for enhancing the interpretability of stochastic programming solutions using multiparametric programming inspired by Rathi et al. [2].

By using the Benders subproblem to generate explicit critical regions, we achieve an objective, model-driven clustering of scenarios. Notably, while the presented case study demonstrated an LP formulation, this framework is fully applicable to problems involving Mixed-Integer Linear Programming (MILP) first-stage decisions.

Applied to a supply chain case study, this method successfully reduced 100 scenarios into 6 logical operating modes. It allowed us to explain the optimal first-stage production levels by rigorously quantifying the trade-off between guaranteed production costs and the expected recourse costs associated with specific critical regions. While this approach is effective, a key limitation is its scalability with respect to problem size. Specifically, as the number of uncertain parameters increases, the multiparametric solution becomes computationally intensive because the total number of critical regions grows combinatorially. Future work will focus on alternative mp algorithms, which generate only the relevant Critical Regions encountered by the Benders cuts, rather than exploring the entire parameter space, enabling scalability to larger networks.

## ACKNOWLEDGEMENTS

This work is based upon work supported by the National Science Foundation under grant no. CMMI-2328160.

Furthermore, we would like to acknowledge the assistance provided by large language models (LLMs), including OpenAI's ChatGPT and Google's Gemini, in the preparation of this manuscript. These tools were used to improve the grammar and fluency of the text, as well as to assist with coding and commenting.

## AUTHOR IDENTIFIERS

Author ORCIDs:
Brahmbhatt P: 0009-0006-1502-1359
Avraamidou S: 0000-0002-9334-9951

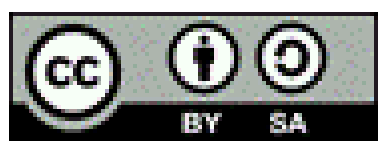